\documentclass[12pt]{article}
\usepackage{a4wide}
\usepackage{mathptmx}
\usepackage{helvet} \usepackage{courier} \usepackage{type1cm}
\usepackage{makeidx} \usepackage{graphicx} \usepackage{multicol}
\usepackage[bottom]{footmisc} \usepackage{subfigure} \usepackage{url}
\usepackage{pstricks} \usepackage{pst-node}
\usepackage{epsfig,hyperref,xspace,todonotes,color}
\newcommand{\bm}[1]{\mbox{$\mathbf #1$}}

\newcommand{\bY}{{\bm Y}} \newcommand{\bM}{{\bm M}}
 
\newcommand{\half}{\frac 1 2} \newcommand{\p}{\mbox{P}\xspace}
\newcommand{\q}{\mbox{Q}\xspace} \newcommand{\pc}{\mbox{PC}\xspace}
\newcommand{\rr}{\mbox{RR}\xspace}
\newcommand{\eqref}[1]{\mbox{(\ref{eq:#1})}}

\newcommand{\ie}{{\em i.e.\/}\xspace}
  \newcommand{\etc}{{\em etc.\/}\xspace}
\newcommand{\etal}{\emph{et al.}\ }

\newcommand{\secref}[1]{\S$\,$\ref{sec:#1}}
\newcommand{\tabref}[1]{Table~\ref{tab:#1}}
\newcommand{\figref}[1]{Figure~\ref{fig:#1}}

\newcommand{\cip}{\mbox{$\,\perp\!\!\!\perp\,$}}
\newcommand{\setto}{\leftarrow}

\makeindex
\begin{document}
\thispagestyle{empty}

\title{Bounding the Probability of Causation in Mediation Analysis}

\author{A.~P.~Dawid\thanks{University of Cambridge {\tt
      apd@statslab.cam.ac.uk}} \and R.~Murtas\thanks{University of
    Cagliari {\tt ro.murtas@gmail.it}} \and M.~Musio\thanks{University
    of Cagliari {\tt mmusio@unica.it}}}


\maketitle

\abstract{Given empirical evidence for the dependence of an outcome
  variable on an exposure variable, we can typically only provide
  bounds for the ``probability of causation'' in the case of an
  individual who has developed the outcome after being exposed.  We
  show how these bounds can be adapted or improved if further
  information becomes available.  In addition to reviewing existing
  work on this topic, we provide a new analysis for the case where a
  mediating variable can be observed.  In particular we show how the
  probability of causation can be bounded when there is no direct
  effect and no confounding.}

\noindent Keywords: Causal inference, Mediation Analysis, Probability
of Causation

\section{Introduction}
\label{sec:introduction}

Many statistical analyses aim at a causal explanation of the data.  In
particular in epidemiology many studies are conducted to try to
understand if and when an exposure will cause a particular disease.
Also in a Court of Law, when we want to assess legal responsibility we
usually refer to causality.  But when discussing this topic it is
important to specify the exact query we want to talk about.  For
example it may be claimed in court that it was Ann's taking the drug
that was the cause of her death.  This type of question relates to the
cause of an observed effect (``CoE'') and is fundamental to the
allocation of responsibility.  On the other hand much of classical
statistical design and analysis, for example randomized agricultural
or medical experiments, has been created to address questions about
the effects of applied causes (``EoC'').  When we address an EoC
query, we are typically asking a hypothetical question: ``What would
happen to Ann if she were to take the drug?''.  At the very same time
we can address alternative hypothetical questions: ``What would happen
to Ann if she were not to take the drug?''.

Assessing the effects of causes can be achieved in straightforward
fashion using a framework based on probabilistic prediction and
statistical decision theory \cite{apd:annrev}.  To formalize the
problem, let $X$ be a binary decision variable denoting whether or not
Ann takes the drug, and $Y$ the response, coded as $1$ if she dies and
$0$ if not.  We denote by $P_1$ [resp.\ $P_0$] the probability
distribution of $Y$ ensuing when $X$ is set to the value $1$ [resp.\
$0$].  The two distributions $P_1$ and $P_0$ are all that is needed to
address EoC-type queries: I can compare these two different
hypothetical distributions for $Y$, decide which one I prefer, and
take the associated decision.

The situation is different for a CoE query, where the drug has already
been taken and the outcome observed.  A natural way to address a CoE
question would be to try to imagine what would have happened to Ann
had she not been taken the drug.  In other words, given the fact that
Ann actually took the drug and died, how likely is it that she would
not have died if she had not taken the drug?  We can not address a CoE
query using only the two distribution $P_1$ and $P_0$.  In fact we can
no longer base our approach purely on the probability distribution of
$Y$ and $X$ conditioned on known facts, since we know the values of
both variables ($Y=1$, $X=1$), and after conditioning on that
knowledge there is no probabilistic uncertainty left to work with.
Nevertheless we want an answer.  This query can be approached by
introducing (for any individual) an associated pair of ``potential
responses'' $\bY := (Y(0), Y(1))$, where $Y(x)$ denotes the value of
the response $Y$ that will be realized when the exposure $X$ is set to
$x$ (which we write as $X\setto x$).  Both potential responses are
regarded as existing, simultaneously, prior to the choice of $X$, the
actual response $Y$ then being determined as $Y=Y(X)$.  However, for
any each individual just one of the potential responses will be
observable.  For example, only $Y(1)$ will be observable if in fact
$X\leftarrow 1$; $Y(0)$ will then be {\em counterfactual\/}, because
it relates to a situation, $X\leftarrow 0$, which is contrary to the
known fact $X\leftarrow 1$.

To address the court's query we use the formulation of
\textit{Probability of Causation}, \pc, as given by Pearl in
\cite{pearl1999} (where it is named \textit{Probability of
  Necessity}).  In terms of the triple $(X_A, Y_A(0), Y_A(1))$, we
define the Probability of Causation in Ann's case
as: 
\begin{equation}
  \label{eq:s}
  \pc_A = \p_A(Y_A(0) = 0 \mid X_A=1, Y_A(1)=1)
\end{equation}
where $\p_A$ denotes the probability distribution over attributes of
Ann.  Knowing that Ann did take the drug ($X_A=1$) and the actual
response was recovery ($Y_A=1$), this is the probability that the
potential response $Y_A(0)$, that would been observed had Ann not
taken the drug, would have been different ($Y_A(0)=0$).  But how are
we to get a purchase on this quantity?

Suppose that a good experimental study tested the same drug taken by
Ann, and produced the data reported in \tabref{aspirin}.
\begin{table}[tb]
  \centering
  \begin{tabular}[t]{lccc}
    & Die & Live & Total\\
    Exposed & 30 & 70 & 100\\
    Unexposed & 12 & 88 & 100
  \end{tabular}
  \centering
  \caption{Deaths in individuals exposed and unexposed to the same drug taken by Ann}
  \label{tab:aspirin}
\end{table}
Since our analysis here is not concerned with purely statistical
variation due to small sample sizes, we take proportions computed from
this table as accurate estimates of the corresponding population
probabilities (see \cite{apd/mm/sef:ba} for issues related to the use
of small-sample data for making causal inferences).  Thus we take
\begin{eqnarray*}
  \Pr(Y = 1 \mid X\setto 1) &=& 0.30\\
  \Pr(Y = 1 \mid X\setto 0) &=& 0.12
\end{eqnarray*}
where we use $\Pr$ to denote population probabilities.

We see that, in the experimental population, individuals exposed to
the drug ($X\setto 1$) were more likely to die than those unexposed
($X\setto 0$), by $28$ percentage points.  So can the court infer that
was Ann's taking the drug that caused her death?  More generally: Is
it correct to use such experimental results, concerning a population,
to say something about a single individual?  This
``Group-to-individual'' (G2i) issue is discussed by Dawid \etal
\cite{fitting_science} in relation to the question ``When can Science
be relied upon to answer factual disputes in litigation?''.  It is
there pointed out that in general we can not obtain a point estimate
for $\pc_A$; but we can provide useful information, in the form of
bounds between which this quantity must lie.

In this paper we show how these bounds can be adapted or improved when
further information is available.  In \secref{simple} we consider the
basic situation where we have information only on exposure and
outcome.  In \secref{cov} we bound the probability of causation when
we have additional information on a pre-treatment covariate.
Section~\ref{sec:pearl} considers the situation in which unobserved
variables confound the exposure-outcome relationship.  Finally in
\secref{med} we introduce new bounds for $\pc$ when a mediating
variable can be observed.  Section~\ref{sec:future} presents some
concluding comments.


%
%
%
%

\section{Starting Point: Simple Analysis}
\label{sec:simple}

In this section we discuss the simple situation in which we have
information, as in \tabref{aspirin}, from a randomized experimental
study.  We need to assume that the fact of Ann's exposure, $X_A$, is
independent of her potential responses $\bY_A$:
\begin{equation}
  \label{eq:sufft}
  X_A \cip \bY_A. 
\end{equation}
Property~\eqref{sufft} parallels the ``no-confounding'' property $X_i
\cip \bY_i$ which holds for individuals $i$ in the experimental study
on account of randomization.  We further suppose that Ann is
exchangeable\ with the individuals in the experiment, \ie she could be
considered as a subject in the experimental population.

On account of \eqref{sufft} and exchangeability, \eqref{s} reduces to
$\pc_A = \Pr(Y(0) = 0 \mid Y(1)=1)$; but we can not fully identify
this from the data.  In fact we can never observe the joint event
$(Y(0)=0;Y(1)=1)$, since at least one of $Y(0), Y(1)$ must be
counterfactual.  In particular, we can never learn anything about the
dependence between $Y(0)$ and $Y(1)$.  However, even without making
any assumptions about this dependence, we can derive the following
inequalities \cite{apd/mm/sef:ba}:
\begin{equation}
  \label{eq:generic}
  1 - \frac 1 \rr \leq
  \pc_A
  \leq
  \frac{\Pr(Y = 0 \mid X\setto 0)}{\Pr(Y=1 \mid X\setto 1)}
\end{equation}
where
$$\rr = \frac{\Pr(Y = 1 \mid X\setto 1)}{\Pr(Y=1 \mid X\setto 0)}$$
is the {\em experimental risk ratio\/} between exposed and unexposed.
And these bounds can be estimated from the experimental data using the
population death rates computed in \secref{introduction}.

In many cases of interest (such as \tabref{aspirin}), we will have
\begin{displaymath}
  \Pr(Y= 1 \mid X \setto 0) <
  \Pr(Y= 1 \mid X \setto 1)
  < \Pr(Y= 0 \mid X \setto 0).
\end{displaymath}
Then the lower bound in \eqref{generic} will be non-trivial, while the
upper bound will exceed 1, and hence be vacuous.

We see from \eqref{generic} that whenever $\rr > 2$ the Probability of
Causation $\pc_A$ will exceed 50\%.  In a civil court this is often
taken as the criterion to assess legal responsibility ``on the balance
of probabilities'' (although the converse is false: it would not be
correct to infer $\pc_A < .5$ from the finding $\rr < 2$).  Since, in
\tabref{aspirin}, the exposed are $2.5$ times as likely to die as the
unexposed ($\rr = 30/12 = 2.5$), we have enough confidence to infer
causality in Ann's case: We have $0.60 \leq \pc_A \leq 1$.

\section{Additional Covariate Information}
\label{sec:cov}
In this Section we show how we can refine the bounds of
\eqref{generic} if further information about a pre-treatment covariate
$S$ is available.  We now take the assumptions of \secref{simple} to
hold after conditioning on $S$ (indeed in cases where the original
assumptions fail, it may well be possible to reinstate them by
conditioning on a suitable covariate $S$). In particular, $X_A \cip
\bY_A \mid S_A$, and $X_i \cip \bY_i \mid S_i$: adjusting for $S$ is
enough to control for confounding, both for Ann and in the study.

\subsection{Fully observable}
\label{sec:fully}
Consider first the situation where we can observe $S$ both in the
experimental data and in Ann.  We can apply the analysis of
\secref{simple}, after conditioning on $S$, to obtain the estimable
lower bound:
\begin{equation}
  \label{eq:covpc}
  1-\frac{1}{\rr(s_A)}  \leq  \pc_A
\end{equation}
where
$$\rr(s) = \frac{\Pr(Y = 1 \mid X\setto 1, S=s)}{\Pr(Y=1 \mid X\setto 0, S=s)},$$
and $s_A$ is Ann's value for $S$.

\subsection{Observable in data only}
\label{sec:dataobs}
But even when we can only observe $S$ in the population, and not in
Ann, we can sometimes refine the bounds in \eqref{generic}.  Thus
suppose $S$ is binary, and from the data we infer the following
probabilities (which are consistent with the data of
\tabref{aspirin}):
\begin{eqnarray}
  \label{eq:pp1}
  \p_A(S=1) &=& 0.50\\
  \label{eq:pp2}
  \p_A(Y=1 \mid X\setto 1, S=1) &=& 0.60\\
  \label{eq:pp3}
  \p_A(Y=1 \mid X\setto 0, S=1) &=& 0\\
  \label{eq:pp4}
  \p_A(Y=1 \mid X\setto 1, S=0) &=& 0\\
  \label{eq:pp5}
  \p_A(Y=1 \mid X\setto 0, S=0) &=& 0.24
\end{eqnarray}

Since we know $X_A = 1, Y_A=1$, from \eqref{pp4} we deduce $S_A=1$,
and so $Y_A(0) = 0$ by \eqref{pp3}.  That is, in this special case we
can infer causation in Ann's case---even though we have not directly
observed her value for $S$.

More generally (see \cite{apd:aberdeen}) we can refine the bounds in
\eqref{generic} as follows:

\begin{equation}
  \label{eq:better}
  \frac{\Delta}{\Pr(Y = 1 \mid X \setto 1)} \leq \pc \leq 1 -
  \frac{\Gamma}{\Pr(Y = 1 \mid X \setto 1)}
\end{equation}
where
\begin{eqnarray*}
  \Delta &=& \sum_s  \Pr(S=s) \times\max\left\{0, \Pr(Y= 1
    \mid X \setto 1, S=s) - \Pr(Y = 1 \mid X \setto 0, S=s)\right\}\\
  \Gamma &=&  \sum_s \Pr(S=s) \times \max\left\{0, \Pr(Y = 1 \mid X \setto 1,
    S=s) - \Pr(Y = 0 \mid X \setto 0, S=s)\right\}
\end{eqnarray*}

These bounds are never wider than those obtained from \eqref{generic},
whicn ignores $S$.

\section{Unobserved Confounding}
\label{sec:pearl}

So far we have assumed no confounding: $X \cip \bY$ (perhaps
conditionally on a stuitable covariate $S$), both for Ann and for the
study data.  Now we drop this assumption for Ann.  Then the
experimental data can not be used, by themselves, to learn about
$\pc_A = \p(Y_A(0) = 0 \mid X_A = 1, Y_A(1)=1)$.

We might however be able to gather additional {\em observational\/}
data, having the same dependence between $X$ and $\bY$ as for Ann.
Let $\q$ denote the joint observational distribution of $(X,Y)$,
estimable from such data.  Tian and Pearl \cite{tianpearl2000} obtain
the following bounds for $\pc_A$, given both experimental and
nonexperimental data:

\begin{equation}
\label{eq:ss}
  {\max\left\{
      0, \frac{\q(Y=1)-\Pr(Y=1 \mid X \setto 0)}{\q(X=1,Y=1)}
    \right\} \leq \pc_A }
  \leq
  \min\left\{
    1, \frac{\Pr(Y=0 \mid X \setto 0)-\q(X=0,Y=0)}{\q(X=1,Y=1)}
  \right\}.
\end{equation}


For example,  suppose that, in addition to the data of \tabref{aspirin},
we have observational data as in \tabref{aspirin2}.

\begin{table}[tb]
  \centering
  \begin{tabular}[t]{lccc}
    & Die & Live & Total\\
    Exposed & 18 & 82 & 100\\
    Unexposed & 24 & 76 & 100
  \end{tabular}
  \caption{Non experimental data}
  \label{tab:aspirin2}
\end{table}

Thus
\begin{eqnarray*}
  \q(Y=1)&=&0.21\\
  \q(X=1,Y=1)&=&0.09\\
  \q(X=0,Y=0)&=&0.38.
\end{eqnarray*}
Also, from \tabref{aspirin} we have $\Pr(Y=1 \mid X \setto 0)=0.12$
(so $\Pr(Y=0 \mid X \setto 0)=1-0.12=0.88$).  From \eqref{ss} we thus
find $1 \leq \pc_A \leq 1$.  We deduce that Ann would definitely have
survived had she not taken the drug.

\section{Mediation Analysis}
\label{sec:med}

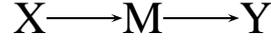
\begin{figure}
  \centering
  \begin{pspicture}(0,-0.5)(3,0.5) \color{black} \rput(0,0){%
      \rnode{1}{\Large X}} \rput(1.5,0){%
      \rnode{2}{\Large M}} \rput(3,0){%
      \rnode{3}{\Large Y}} \ncline[linecolor=black]{->}{1}{2}
    \ncline[linecolor=black]{->}{2}{3}
  \end{pspicture}
  \caption{Directed Acyclic Graph representing a mediator $M$,
    responding to exposure $X$ and affecting response $Y$.  There is
    no ``direct effect'', unmediated by $M$, of $X$ on $Y$.}
  \label{fig:nodirect}
\end{figure}

In this Section we bound the Probability of Causation for a case where
a third variable, $M$, is involved in the causal pathway between the
exposure $X$ and the outcome $Y$.  Such a variable is called a
\textit{mediator}.  In general, the total causal effect of $X$ on $Y$
can be split into two different effects: One mediated by $M$ (the {\em
  indirect effect\/}) and one not so mediated (the {\em direct
  effect\/}).  Here we shall only consider the case of no direct
effect, as intuitively described by \figref{nodirect}.  We shall be
interested in the case that $M$ is observable in the experimental
data, but is not observed for Ann, and see how this additional
experimental evidence can be used to refine the bounds on $\pc_A$.

To formalize our assumption of ``no direct effect'', we introduce
$M(x)$, the potential value of $M$ for $X\setto x$, and $Y(m)$, the
potential value of $Y$ for $M\setto m$ and $X\setto x$, where $x$ is
any value---the irrelevance of that value representing the property
that $X$ has no effect on $Y$ over and above that transmitted through
its influence on the mediator $M$.  The potential value of $Y$ for
$X\setto x$ (in cases where there is no intervention on $M$, which we
here assume) is then $Y^*(x) := Y\{M(x)\}$.

In the sequel we restrict to the case that all variables are binary,
and define $\bM := (M(0), M(1))$, $\bY := (Y(0),Y(1))$, $\bY^*
:=(Y^*(0),Y^*(1))$.  In particular, we have observable variables
$(X,M,Y) = (X, M(X), Y^*(X))$.  We denote the bivariate distributions
of the potential response pairs by
\begin{eqnarray*}
  m_{ab} &:=& \Pr(M(0)=a, M(1)=b)\\
  y_{rs} &:=& \Pr(Y(0)=r, Y(1)=s)\\
  y^*_{rs} &:=& \Pr(Y^*(0)=r,Y^*(1)=s).
\end{eqnarray*}
Then 
\begin{eqnarray*}
  m_{a+} &=& \Pr(M=a\mid X\setto 0)\\
  m_{+b} &=& \Pr(M=b \mid X\setto 1)\\
  y_{r+} &=& \Pr(Y=r \mid M\setto 0)\\
  y_{+s} &=& \Pr(Y=s \mid M\setto 1)\\
  y^*_{r+} &=& \Pr(Y=r \mid X\setto 0)\\
  y^*_{+s} &=& \Pr(Y=s \mid X\setto 1)
\end{eqnarray*}
where $y_{r+}$ denote $\sum_{s=0}^1 y_{rs}$, \etc

In addition to the assumptions of \secref{simple} we further suppose
that none of the causal mechanisms depicted in \figref{nodirect} are
confounded---expressed mathematically by assuming mutual independence
between $X$, $\bM$ and $\bY$ (both for experimental individuals, and
for Ann).
Then $m_{a+}$, $m_{+b}$, $y_{r+}$, $y_{+s}$, $y^*_{r+}$, $y^*_{+s}$
are all estimable from experimental data in which $X$ is randomized,
and $M$ and $Y$ are observed.

It is also then easy to show the Markov property:
$$Y \cip X \mid M.$$
This observable property can serve as a test of the validity of our
conditions.  It implies
\begin{eqnarray}
  \label{eq:*00}
  y^*_{00} &=& m_{00}y_{0+} + (m_{01}+m_{10})y_{00} + m_{11} y_{+0}\\
  \label{eq:*01}
  y^*_{01} &=& m_{01}y_{01} + m_{10}y_{10}\\
  \label{eq:*10}
  y^*_{10} &=& m_{01}y_{10} + m_{10}y_{01}\\
  \label{eq:*11}
  y^*_{11} &=& m_{00}y_{1+} + (m_{01}+m_{10})y_{11}+ m_{11} y_{+1},
\end{eqnarray}
and
\begin{eqnarray}
\label{eq:ystar1}
  y^*_{r+} &=& m_{0+}y_{r+} +m_{1+}y_{+r}\\
\label{eq:ystar2}
  y^*_{+s} &=&  m_{+0}y_{s+} +  m_{+1}y_{+s}
\end{eqnarray}

Suppose now that we observe $X_A=1$, $Y_A=1$, but do not observe
$M_A$.  We have
\begin{equation}
\label{eq:pc}
\pc_A = \frac{y^*_{01}}{y^*_{+1}}=\frac{m_{01}y_{01} + m_{10}y_{10}}{y^*_{+1}}.
\end{equation}

The denominator of \eqref{pc} is $\Pr(Y=1 \mid X\setto 1)$, which is
estimable from the data.

As for the numerator, this can be expressed as
\begin{equation}
   \label{eq:num}
   2\mu\eta + A\mu + B \eta + AB
   = 2(\mu+B/2)(\eta+A/2) + AB/2
\end{equation}
with $\mu = m_{01}$, $\eta= y_{01}$, $A = y_{+0}-y_{0+}$, $B =
m_{+0}-m_{0+}$.  Note that $A$, $B$ are identified from the data,
while for $\mu$ and $\eta$ we can only obtain inequalities:
\begin{displaymath}
  \begin{array}[c]{rcccl}
    \max\{0,-B\} &\leq& \mu &\leq& \min\{m_{0+},m_{+1}\}\\
    \max\{0,-A\} &\leq& \eta &\leq& \min\{y_{0+},y_{+1}\},
  \end{array}
\end{displaymath}
so that
\begin{equation}
  \label{eq:bounds}
  \begin{array}[c]{rcccl}
    |B/2| &\leq& \mu+B/2 & \leq & \min\{\half(m_{0+} + m_{+0}), \half(m_{1+} + m_{+1})\}\\
    |A/2| &\leq& \eta+A/2& \leq & \min\{\half(y_{0+} + y_{+0}), \half(y_{1+} + y_{+1})\}.  
  \end{array}
\end{equation}
The lower [resp., upper] limit for \eqref{num} will be when $\mu+B/2$
and $\eta+A/2$ are both at their lower [resp., upper] limits.  In
particular, the lower limit for \eqref{num} is $\max\{0,AB\}$.  Using
\eqref{ystar1} and \eqref{ystar2}, we compute $AB = y^*_{+1} -
y^*_{1+}$, which leads to the lower bound
\begin{equation}
  \label{eq:lb}
  \pc_A \geq 1-\frac{\Pr(Y = 1 \mid X\setto 1)}{\Pr(Y=1 \mid X\setto 0)}=
  1- \frac{1}{\rr},
\end{equation}
exactly as for the case that $M$ was not observed.  Thus the
possibility to observe a mediating variable in the experimental data
has not improved our ability to lower bound $\pc_A$.

We do however obtain an improved upper bound.  Taking into account the
various possible choices for the upper bounds in \eqref{bounds}, the
uper bound for the numerator of \eqref{pc}, in terms of experimentally
estimable quantities, is given in \tabref{UB}.

\begin{table}[tb]
  \centering
  \begin{tabular}{l|cc}
    & $m_{1+} + m_{+1} \geq 1$ & $m_{1+} + m_{+1} < 1$\\ 
    \hline
    $y_{1+} + y_{+1} \geq 1$ &  $m_{0+}y_{0+} + m_{+0}y_{+0}$ & $m_{1+}y_{+0} + m_{+1}y_{0+}$\\
    $y_{1+} + y_{+1}  <   1$ &  $m_{0+}y_{+1} + m_{+0}y_{1+}$ & $m_{1+}y_{1+} + m_{+1}y_{+1}$
  \end{tabular}
  \caption{Upper bound for numerator of \eqref{pc}}
  \label{tab:UB}
\end{table}

\subsection{Example}
\label{sec:exmed}
Suppose we obtain the following values from the data:
\begin{eqnarray*}
  \Pr(M=1 \mid X\setto 1) &=& 0.25  \\
  \Pr(M=1 \mid X\setto 0) &=& 0.025 \\
  \Pr(Y=1 \mid M\setto 1) &=& 0.9   \\
  \Pr(Y=1 \mid M\setto 0) &=& 0.1
\end{eqnarray*}
(these are consistent with \tabref{aspirin}).

Then we find $0.60 \leq\pc_A \leq 0.76$; whereas without taking
account of the mediator we would have no non-trivial upper bound.

\section{Discussion}
\label{sec:future}

In this paper we have considered estimation of the Probability of
Causation in a number of contexts, including a novel analysis for the
case of a mediating variable, in the absence of a direct effect.  As
we saw in \secref{med} considering such a third variable in the
pathway between exposure and outcome can lead to an improved upper
bound, although conclusions about the lower bound remains the same.

Even if the case of no direct effect is special and unusual, there
certainly do exist cases, such as the relationship between anxiolytics
and cars crash mediated by alcohol consumption, or the relationship
between aspirin and yellow jaundice due to hemolytic anemia mediated
by favism, where this assumption is plausible.  The next step will be
generalize our analysis to more general cases of mediation, allowing
for a direct effect and for unobserved confounding.


\begin{thebibliography}{99}

\bibitem{apd:aberdeen} Dawid,~A.~P. (2011).  The role of scientific
  and statistical evidence in assessing causality.  In {\em
    Perspectives on Causation\/} (R.~Goldberg, Ed.).  Oxford: Hart
  Publishing, 133--147.

\bibitem{apd:annrev} Dawid,~A.~P. (2014).  Statistical causality from
  a decision-theoretic perspective.  {\em Annual Review of Statistics
    and Its Application\/} {\bf 2}.  In Press.

\bibitem{apd/mm/sef:ba} Dawid,~A.~P., Musio,~M., and
  Fienberg,~S.~E. (2014).  From statistical evidence to evidence of
  causality.  \href{http://arxiv.org/abs/1311.7513}{\tt
    arXiv:1311.7513}.

\bibitem{tianpearl2000} Tian,~J. and Pearl,~J. (2000).  Probabilities
  of causation: {B}ounds and identification.  {\em Annals of
    Mathematics and Artificial Intelligence\/} {\bf 28}, 287--313.

\bibitem{pearl1999} Pearl,~J. (1999).  Probabilities of causation:
  {T}hree counterfactual interpretations and identification.  {\em
    Synthese \/} {\bf 121}, 93-149.
 
\bibitem{fitting_science} Dawid,~A.~P., Fienberg,~S. and
  Faigman,~D. (2014a). Fitting science into legal contexts: Assessing
  effects of causes or causes of effects?  {\em Sociological Methods
    and Research\/} {\bf 43}, 359--390.

\end{thebibliography}
\end{document}